\documentclass[12pt]{article}
\usepackage{hhline}
\usepackage{float}
\usepackage{bm}
\usepackage{color}
\usepackage[normalem]{ulem}

\usepackage{latexsym,amsmath, amscd, amssymb} 
\usepackage{subfig}
\usepackage[dvips]{epsfig,graphicx}

\newcommand{\appsection}[1]{\let\oldthesection\thesection
 \renewcommand{\thesection}{Appendix \oldthesection}
 \section{#1}\let\thesection\oldthesection}

\newcommand{\bea}{\begin{eqnarray}}
\newcommand{\eea}{\end{eqnarray}}

\title{Rational Construction of Stochastic Numerical Methods for Molecular Sampling}

\author{Benedict Leimkuhler and Charles Matthews \thanks{School of Mathematics and Maxwell Institute of Mathematical Sciences, James Clerk Maxwell Building, Kings Buildings, University of Edinburgh, Edinburgh, EH9 3JZ, UK}}

\begin{document}
\setlength{\baselineskip} {20pt}


\maketitle

\begin{abstract}
In this article, we focus on the sampling of the configurational Gibbs-Boltzmann distribution, that is, the calculation of averages of functions of the position coordinates of a molecular $N$-body system modelled at constant temperature.    We show how a formal series expansion of the invariant measure of a Langevin dynamics numerical method can be obtained in a straightforward way using the Baker-Campbell-Hausdorff lemma.   We then compare Langevin dynamics integrators in terms of their invariant distributions and demonstrate a superconvergence property (4th order accuracy where only 2nd order would be expected) of one method in the high friction limit; this method, moreover, can be reduced to a simple modification of the Euler-Maruyama method for Brownian dynamics involving a non-Markovian (coloured noise) random process.    In the Brownian dynamics case, 2nd order accuracy of the invariant density is achieved. All methods considered are efficient for molecular applications (requiring one force evaluation per timestep) and of a simple form.   In fully resolved (long run) molecular dynamics simulations, for our favoured method,  we observe up to two orders of magnitude improvement in configurational sampling accuracy for given stepsize with no evident reduction in the size of the largest usable timestep compared to common alternative methods. 

\end{abstract}

\par\noindent {\bf keywords:} molecular dynamics; sampling; Langevin dynamics; Brownian dynamics; stochastic dynamics.

\newpage 
\section{Introduction}

Let $U:{\bf R}^N\rightarrow {\bf R}$ be the potential energy function of a classical model for a molecular system.   A fundamental challenge  is to sample the configurational Gibbs-Boltzmann (canonical) distribution with density
\begin{equation}
\bar{\rho}_{\beta}(x) = Z^{-1} \exp(-\beta U(x)), \label{configdist} 
\end{equation}
where $\beta^{-1} = k_BT$ where $k_B$ is Boltzmann's constant, $T$ is temperature, and $Z$ is  a normalization constant so that $\bar{\rho}_{\beta}$ has unit integral over the entire configuration space.  A wide variety of methods are available to calculate averages with respect to $\bar{\rho}_{\beta}$; among these, some of the most popular are based on Brownian dynamics or Langevin dynamics (defined in the  phase space of positions and momenta) \cite{BBK, Ta1995, MiSc96, MiTr2004, BuLy2009}. { (In this article we focus exclusively on molecular dynamics techniques; molecular models can also be sampled using Monte-Carlo methods, and, more generally, using hybrid algorithms which combine molecular dynamics with a Metropolis-Hastings test in order to correct averages. For a recent review of such schemes see \cite{FreeEnergy}.)}  Recall that Brownian dynamics (overdamped Langevin dynamics) is a system of It\={o}-type stochastic differential equations of the form 
\begin{equation} \label{BD}
{\rm d} x = - M^{-1} \nabla U(x) {\rm d} t + \sqrt{2\beta^{-1}}M^{-1/2}{\rm d} W, \hspace{0.2in} x(0)=x_0,
\end{equation}
where ${\rm d}W(t)$ is the infinitesimal increment of a vector of stochastic Wiener processes $W(t)$, and  $M$ is a positive  (we assume here diagonal) mass matrix.\footnote[1]{In terms of sampling the Gibbs distribution, $M$ is in fact arbitrary, but it may be useful to allow for coordinate scaling.}
A simple and popular method for numerical solution of Eq. \ref{BD} is the Euler-Maruyama method
\[
x_{n+1} = x_n - hM^{-1}\nabla U(x_n)  + \sqrt{2k_BT h}M^{-1/2}R_n,
\]
where $R_n$ is a vector of random variables with standard normal distribution.   This produces a sequence of points $x_0, x_1, x_2, \ldots$,  which, following a certain relaxation period, are approximately distributed according to the canonical invariant distribution.   Euler-Maruyama  has the property that the time averages along discrete trajectories, in the limit of large time ({under appropriate conditions on the potential $U(x)$} and assuming no effects from floating point rounding error), have error proportional to $h$.      One of the observations of this article is that the simple modification 
\begin{equation}
x_{n+1} = x_n - h M^{-1} \nabla U(x_n) + \sqrt{\frac{k_BT h}{2}} M^{-1/2}(R_n+R_{n+1}) \label{BD-mod}
\end{equation}
provides a second order approximation of stationary averages.     

We arrive at this scheme by considering the large friction limit in a particular numerical method for Langevin dynamics.     Recall that Langevin dynamics
is a stochastic-dynamical system involving both positions and momenta $p$ of the form
\begin{equation}
{\rm d} x =  M^{-1} p{\rm d}t, \hspace{0.1in} {\rm d} p =  [-\nabla U(x) -\gamma p]{\rm d}t + \sigma M^{1/2}{\rm d}W, \label{lang}
\end{equation}
where $W=W(t)$ is again a vector of $N$ independent Wiener processes,  and $\gamma>0$ is a free parameter, the friction coefficient. 
The methods that are in fact the primary focus of this paper are splitting integrators that decompose the stochastic vector field of Langevin dynamics into simpler vector fields which can be solved exactly.  The composition method that results cannot be directly related to a stochastic differential equation, so the analogy with backward error analysis for deterministic problems \cite{HaLuWa2006,LeRe2005} is incomplete, but nonetheless the invariant measure associated to the numerical method can be derived, using the Baker-Campbell-Hausdorff (BCH) expansion, as an asymptotic series in two parameters: the stepsize and the reciprocal of the friction coefficient.   The superconvergence result alluded to above is then obtained in the high friction limit.   

The Langevin stepsize must be understood to be proportional to the square root of the stepsize that appears in Eq.\ \ref{BD-mod}, so in Langevin dynamics an effective 4th order approximation is obtained,  but only for the marginal configurational invariant distribution, { Eq.} \ref{configdist}.  Our approach also provides a simple method for comparative assessment of the invariant measure of a class of Langevin integrators.    

Molecular dynamics is a large family of modelling techniques which is widely used in different application areas and for different purposes \cite{FrSm2001}. This article is addressed specifically to the topic of calculating averages with respect to an invariant distribution, and will probably be of highest interest for applications in molecular sampling, in particular the calculation of averages with respect to the configurational density Eq. \ref{configdist}.  The high friction limit renders Langevin dynamics unsuitable for dynamical modelling (except as a method for generating starting configurations for dynamical exploration).   It is worth noting that invariant measure computations arise frequently in applications other than molecular modelling, and the techniques described here would be of potential use in many of these.

{In large scale simulations, it should be understood that the statistical error (dependent on the number of samples used) is typically  the dominant concern.  Our approach focuses on the truncation error of the invariant distribution, thus the greatest benefit would be seen only when the statistical error is well controlled.  Nonetheless we observe in our numerical experiments that the least biased scheme from the point of view of the error introduced in configurational sampling is also as efficient as the alternatives and the most robust with respect to variation of the parameters (stepsize, friction coefficient), thus there is effectively no price for the improvement in accuracy.   Moreover, it  would be possible to complement the methods proposed here by  procedures such as importance sampling \cite{ImportanceSampling} to further reduce the statistical error.}

With regard to the approximation of canonical averages, methods have previously been constructed for Brownian dynamics with order $>1$ and for Langevin dynamics with order $>2$ \cite{Ta1995,MiTr2004}, but these require multiple evaluations of the force; for this reason they are not normally viewed as competitive alternatives for molecular sampling  \cite{LaMaLe2007}.     By contrast all of the methods described in this article use a single force evaluation at each timestep.

The approach used here may be compared to other recent works on stochastic numerical methods, and, in particular  \cite{SkIz02,MiTr2003,Sh2006,Me2007, BoOw2010,Da2010,Zy2011}.  Our technique differs from these in (i) the direct focus on the stationary configurational distribution, and (ii) the use of the BCH expansion.    Other articles (see e.g. \cite{BuLy2009}) which address the invariant measure in Langevin-type stochastic differential equations do not use the backward error analysis (and do not find the superconvergent scheme).  

{ The rest of this article is as follows. In Section 2 we review necessary background on stochastic differential equations for sampling from the Gibbs measure. Section 3 presents our expansion of the associated perturbed invariant measure and calculations involving the use of the Baker-Campbell-Hausdorff theorem. Section 4 describes the reduction of the methods in the case of overdamped Langevin Dynamics. Section 5 demonstrates the theory obtained using numerical experiments to verify the results.}

\section{Background}
The Ornstein-Uhlenbeck (OU) stochastic differential equation is an It\={o} equation of the form
\[
{\rm d}u = -\gamma u {\rm d}t + \sigma {\rm d} W,
\] 
where $u$ is a random variable defined for each time $t$, $\gamma$ and $\sigma$ are positive parameters and ${\rm d}W$ represents the infinitesimal increment of a Wiener process.      Langevin dynamics (Eq. \ref{lang}) combines the Ornstein-Uhlenbeck stochastic vector field with conservative dynamics.

For stochastic differential equations, the density evolves according to an evolution equation (the Fokker-Planck or forward Kolmogorov equation) of the form
\[ 
\frac{\partial \rho}{\partial t} = {\cal L}^* \rho, 
\]
where ${\cal L}^*$ is a second order differential operator.   In the case of Langevin dynamics, the relevant Kolmogorov operator is defined by its action on a function $\phi$ of the variables of the system by
\[
{ {\cal L}_{\rm LD}^*\phi = -M^{-1}p \cdot \nabla_x \phi  + \nabla U(x) \cdot \nabla_p \phi+ \gamma  \nabla_p\cdot (p \phi) + \frac{\sigma^2}{2}\Delta_p \phi}, 
\]
where $\Delta_p$ is the mass-weighted Laplacian in the momenta: $\Delta_p = \sum_i m_i \frac{\partial ^2}{\partial p_i^2}$. 
By choosing $\sigma = \sqrt{2k_B T\gamma}$ one easily checks that the {Gibbs distribution with density $\rho_{\beta}=\tilde{Z}^{-1} \exp(-\beta H)$ is a steady state of the Kolmogorov equation, where $\tilde{Z}$ is a normalising constant} and $H(x,p) = p^TM^{-1}p/2 + U(x)$ represents the Hamiltonian energy function.
For later reference, we define averages with respect to the Gibbs distribution by
\[
 { \mathbb{E} \left( \phi(x,p)  \right) = \int \phi(x,p) \rho_\beta  (x,p) \, {\rm d} x \,  {\rm d}p }.
\]
Assuming $U$ is $C^{\infty}$ it is possible to demonstrate that the operator ${\cal L}_{\rm LD}^*$ is hypoelliptic by using  H\"ormander's criterion \cite{Ho} based on iterated commutators, and this implies that the Gibbs measure is the unique steady state (up to normalization).  Even stronger is the result of \cite{MaStHi2002} which demonstrates the existence of a spectral gap for the operator ${\cal L}_{\rm LD}$.   Many of the challenges related to obtaining formal analytical results for stochastic differential equations relate to singularities of the potential and/or the assumption of an unbounded solution domain.  However, with periodic boundary conditions and strong repulsive potentials (e.g. Lennard-Jones potentials) we observe that configurations typically evolve in a bounded set and remain far from singular points (the radial distribution vanishes in a large interval around the origin).  Indeed it is a simple calculation to demonstrate that for a Lennard-Jones system with potentials $\varphi(r) = (\sigma/r)^12-2(\sigma/r)^6$ the expected number of samples required to observe $r \in (0,\sigma/2)$ at unit temperature would   involve a simulation of duration far greater than the age of the universe, due to the steepness of the potential close to the origin. In our simulations of a small Lennard-Jones cluster in Section 5, we did not observe a separation of two atoms beneath $0.7\sigma$ in all of the nearly $10^{11}$ timesteps performed to gather statistics; a separation less than $0.5\sigma$ could only be seen at very large stepsizes, when instabilities due to other components of model, e.g. harmonic bonds, would have anyway rendered a typical molecular simulation useless.  Thus it is somewhat of a moot point whether we simply assume that configurations stay well away from singular points (domain restriction) or that the potential has been smoothly cut off to remove the singularity; in no case will we encounter the atomic collision singularity in a simulation of the type envisioned here.

For the purpose of deriving practical methods we assume that (i) the positions are confined to a periodic simulation box $\Omega=L_x \mathbb{T} \times L_y \mathbb{T} \times L_z \mathbb{T}$, where \mbox{$\mathbb{T}=\mathbb{R}/\mathbb{Z}$} is the one-dimensional torus, and (ii) $U$ is $C^{\infty}$ on $\Omega$.  These assumptions, which are realistic in most molecular dynamics applications for the reasons mentioned above, allow us to use recent results in hypocoercivity  \cite{Vil2009, PavHai2008, JouSto2012} to establish the regularity properties of ${\cal L}_{\rm LD}$.   Specifically, we have as a consequence of these articles, that
\begin{itemize}
\item ${\cal L}_{\rm LD}^* f =0$ has a unique solution  in $C^{\infty}(\Omega)$, i.e., the Gibbs density $\rho_{\beta}$,
\item ${\cal L}_{\rm LD}^*$ has a compact resolvent \cite{Vil2009, PavHai2008}, and the Gibbs state is therefore exponentially attracting.
\end{itemize}
Note that, if, for some given $g$,  there are two solutions of ${\cal L }_{\rm LD}^*f =g$ then, using the linearity of the operator 
 these may differ only in a constant scaling of the Gibbs density.

%



\subsection{Timestepping Methods.}
In a splitting method for a deterministic system $\dot{z} = f(z)$, one divides the vector field $f$ into exactly solvable parts, i.e. $f=f_1+f_2$, which are treated sequentially within a timestep.    An example of such a splitting method for the Hamiltonian system with energy $H(x,p) =\sum_i p_i^2/(2m_i) + U(x_1,x_2,\ldots,x_N)$ is the ``symplectic Euler'' method defined by $f_1 = \sum_i m_i^{-1}p_i \partial_{x_i}$, $f_2 = -\sum_i \frac{\partial U}{\partial x_i} \partial_{p_i}$.
By dividing the vector field as $f=\frac{1}{2} f_1 + f_2 + \frac{1}{2} f_1$, solving each vector field in turn, we obtain the so-called position Verlet method, and by switching the roles of $f_1$ and $f_2$ we obtain the velocity Verlet method.  Splitting methods like these are explicit and this feature is of particular importance in molecular dynamics, where the force calculation is the usual measure of per-timestep computational complexity.   

In a similar way, Langevin dynamics may be treated by splitting \cite{SkIz02,Me2007,BoOw2010}.  
For example, one may divide the Langevin system (Eq.\ \ref{lang}) into three parts
\begin{equation}
{\rm d} \left [ \!\!\begin{array}{c} x\\p\end{array}\!\! \right ]=
\underbrace{\left [ \!\!\begin{array}{c} M^{-1} p \\0\end{array}\!\! \right ]{\rm d}t}_{\rm A} +
\underbrace{\left [\!\! \begin{array}{c} 0 \\-\nabla U\end{array}\!\! \right ]{\rm d}t}_{\rm B} +
\underbrace{\left [\!\! \begin{array}{c} 0\\ -\gamma p{\rm d}t +\sigma M^{1/2}{\rm d} W \end{array} \!\!\right ],}_{\rm O}
\label{gla2}
\end{equation}
and each of the three parts may be solved `exactly'.  In the case of the OU part (labelled here simply as O) we mean by this that we realize the stochastic process by the equivalent formula
\begin{equation}\label{OUsolution}
p(t) = e^{-\gamma t} p(0) + \sigma \sqrt{ (1- e^{-2\gamma t}) \beta^{-1} } M^{1/2} R(t),
\end{equation}
where $R(t)$ is a vector of uncorrelated independent standard normal random processes (white noise).      One method based on the splitting in Eq.\ \ref{gla2} is defined by the composition
\[
\psi_{\rm ABAO}^{\delta t} = \exp((\delta t/2){\cal L}_{\rm A}) \exp(\delta t{\cal L}_{\rm B}) \exp((\delta t/2){\cal L}_{\rm A})
\exp(\delta t {\cal L}_{\rm O}),
\]
where $\exp(\delta t {\cal L}_f)$ represents the phase space propagator associated to the (deterministic or stochastic) vector field $f$, and we use $\delta t$ as the timestep in Langevin dynamics (later we use $h$ for the timestep of an associated Brownian dynamics).   The deterministic part is approximated by the position Verlet method.   This is referred to as a Geometric Langevin Algorithm of order two (GLA-2) following \cite{BoOw2010}; an alternative is to use velocity Verlet for the Hamiltonian part ($\psi_{\rm BABO}^{\delta t}$).  
A simple generalization of GLA methods is obtained by interspersing integrators associated to parts of the Hamiltonian vector field with exact OU solves, thus we have 
\[
\psi_{\rm ABOBA}^{\delta t} = \exp((\delta t/2){\cal L}_{\rm A}) \exp((\delta t/2){\cal L}_{\rm B}) \exp(\delta t {\cal L}_{\rm O}) \exp((\delta t/2){\cal L}_{\rm B}) \exp((\delta t/2){\cal L}_{\rm A}),
\]
and $\psi_{\rm BAOAB}$ defined in an analogous way. Recent work in \cite{SiChCr} similarly uses exact OU solves to give an integrator equivalent to $\psi_{\rm OBABO}$ (in our notation), though with a reparameterised timestep. The analysis technique we use can employed to study  many such integrators, though for brevity we will limit this article to discussion only on a select few interesting cases.

An alternative  integrator termed the Stochastic Position Verlet (SPV) method \cite{SkIz02,Me2007},
relies on the splitting
\[
{\rm d} \left [ \begin{array}{c} x\\p\end{array} \right ] = \left [ \begin{array}{c} M^{-1} p \\0\end{array} \right ]{\rm d}t +
\left [ \begin{array}{c} 0 \\-\nabla U{\rm d}t -\gamma p{\rm d}t +\sigma M^{1/2} {\rm d} w \end{array} \right ].
\]
SPV is not a (generalized) GLA-type method, although it, as each of  the generalized GLA schemes, is {\em quasisymplectic} in the language of \cite{MiTr2003}.   Likewise the commonly used method of Brunger, Brooks and Karplus (BBK) \cite{BBK} is not of the (generalized) GLA family.   Details of all methods examined are given in the Appendix.

To slightly simplify the presentation that follows, we make the change of variables $q \rightarrow M^{-1/2}q$, $p\rightarrow M^{+1/2}p$, {with a corresponding adjustment of the potential; this is equivalent to assuming $M=I$.}

\section{Expansion of the Invariant Measure}
{We shall work here with formal series expansions, however we expect that our derivation could be rigorously founded using techniques found in \cite{Tal2002} and \cite{Mat2010}.}  Associated to any given splitting-based method for Langevin dynamics, we define the operator $\hat{\cal L}^*$ that characterises the propagation of density by an expansion of the form:
\begin{equation} \label{opexp}
\hat{\cal L}^* = {\cal L}_{\rm LD}^* + \delta t {\cal L}_1^* + \delta t^2 {\cal L}_2^* + O(\delta t^3),
\end{equation}
For example, for the method labelled by the string ABAO, we have
\[
\exp(\delta t \hat{\cal L}^*_{\rm ABAO})  = \exp((\delta t/2) {\cal L}^*_{\rm A}) \exp(\delta t {\cal L}^*_{\rm B}) \exp((\delta t/2) {\cal L}^*_{\rm A}) \exp(\delta t {\cal L}^*_{\rm O}).
\]
The perturbation series may be found by successive applications of the BCH  expansion \cite{HaLuWa2006} and linearity properties of the Kolmogorov operator. {However, unlike in the deterministic case, the terms that appear in the series cannot be associated to modified vector fields or even SDEs \cite{DeFa2011}.}

Note that, when iterated $n+1$ times, the method ABAO produces a sequence of the form
\[
e^{ \frac{\delta t}{2} {\cal L}_{\rm A}}e^{\frac{\delta t}{2} {\cal L}_{\rm B}}[ \psi_{\rm BAOAB}^{\delta t}]^n e^{\frac{\delta t}{2} {\cal L}_{\rm B}}e^{\frac{\delta t}{2} {\cal L}_{\rm A}} e^{ {\delta t} {\cal L}_{\rm O} }
\]
thus, with a minor coordinate transformation, the dynamics sample the same invariant density as BAOAB.   Similarly BABO and OBAB are essentially the same method as ABOBA.    For this reason we concentrate in the remainder of this article on ABOBA and BAOAB.  For these two methods, the symmetry implies that the odd order terms in Eq. \ref{opexp} vanish identically using the Jacobi identity in the BCH expansions.

After deriving $\hat{\cal L}^*$ in this way, we seek the invariant distribution which satisfies $\hat{ {\cal L} }^* \hat{\rho} = 0$. For the BAOAB and ABOBA methods, we make the ansatz that the invariant measure of the numerical method has the simple form
\begin{equation} \label{ansatz}
\hat{\rho} \propto \exp(-\beta [ H + \delta t^2 f_2 + \delta t^4 f_4 + \ldots]).
\end{equation}
Although some technical issues might be encountered, we believe that the existence of such an expansion can be made rigorous using techniques found in \cite{JouSto2012}, based on the regularity of the operator ${\cal L}^*_{\rm LD}$.     We may rewrite this as 
\[
\hat{\rho} \propto \rho_\beta\left(1 -\beta \delta t^2 f_2(x,p) + O \left(\delta t^4\right) \right).
\]
This means that the equation $\hat{\cal L}^* \hat{\rho} =0$  becomes
\[
\left( {\cal L}_{\rm LD}^* + \delta t^2 {\cal L}_2^* + \ldots \right) \left( \rho_\beta - \delta t^2 \beta \rho_\beta f_2 + \ldots \right) = 0.
\]
Equating second order terms in $\delta t$ gives
\begin{equation}
{\cal L}_{\rm LD}^* \left( \rho_\beta f_2 \right) = \beta^{-1} {\cal L}_{2}^* \rho_\beta . \label{eqn::pde}
\end{equation}

The equation ${\cal L}_{\rm LD}^* g = 0$ has unique solution $g=\rho_\beta$, up to a constant multiple. Hence the homogeneous solution to the above PDE is $f_2(x,p) = c$, for some constant $c$; we therefore require a particular solution $f_2$ of Eq. \ref{eqn::pde}. According to the Fredholm Alternative,   the equation has a solution provided that, for any solution of
\[
{\cal L}_{\rm LD} g = 0,
\]
we have
\[
\int g {\cal L}_{2}^* \rho_\beta  =0.
\]
As the only solutions of ${\cal L}_{\rm LD}g=0$ are the constants, we require
\begin{equation} \label{FA}
\int {\cal L}_{2}^* \rho_\beta = 0.
\end{equation}

\subsection{Calculation of the inhomogeneity}
For a symmetric splitting method such as the BAOAB method, recall that we can use the Baker-Campbell-Hausdorff \cite{HaLuWa2006}  formula to find the overall one-step perturbation operator for the scheme. For linear operators $X$, $Y$ and $Z$, we have the relation
\[
e^{\frac{\delta t}{2} X } e^{\frac{\delta t}{2} Y } e^{\delta t Z} e^{\frac{\delta t}{2} Y} e^{\frac{\delta t}{2} X } = e^{\delta t S},
\]
where
\begin{align*}
S =&  X + Y + Z + \frac{\delta t^2}{12} \left(\left[Z, \left[ Z, Y + X \right] \right] \right. \left. + \left[Y, \left[ Y , X \right] \right] + \left[Z, \left[ Y , X \right] \right]+ \left[Y, \left[ Z , X \right] \right]  \right. \\
& \left. - \tfrac{1}{2} \left[Y, \left[ Y, Z \right] \right] - \tfrac{1}{2} \left[X, \left[ X, Z \right] \right] - \tfrac{1}{2} \left[X, \left[ X, Y \right] \right]  \right) + {\cal O}(\delta t^4).
\end{align*}
Here $[X,Y] = XY - YX$ is the commutator of $Y$ and $X$.

In the case of the BAOAB method, we take
\[
X = {\cal L}^*_{\textrm B} = \nabla U(x) \cdot \nabla_p , \, \, \, Y = {\cal L}^*_{\textrm A} = -p \cdot \nabla_x ,
\]
and
\[
Z  = {\cal L}^*_{\textrm O} =  \gamma  \nabla_p\cdot (p \cdot) + \frac{\sigma^2}{2}\Delta_p,
\]
to compute the perturbed operator for the method. A similar analysis can be conducted for the other generalised GLA-type methods considered in this paper.

To compute $\hat{\cal L}^*$ we simply plug in our choices into the BCH formula to obtain:
\begin{eqnarray*}
\hat{{\cal L}^* } &=& {\cal L}^*_{\rm LD} + \delta t^2 {\cal L}^*_2 + O \left( \delta t^4 \right), \\
&=& {\cal L}^*_{\rm LD} + \frac{\delta t^2}{12} \left( \left[{\cal L}^*_{\rm O}, \left[ {\cal L}^*_{\rm O}, {\cal L}^*_{\rm Det} \right] \right]  + \left[{\cal L}^*_{\rm A}, \left[ {\cal L}^*_{\rm A} , {\cal L}^*_{\rm B} \right] \right] \right. \\
&& \left. + \left[{\cal L}^*_{\rm O}, \left[ {\cal L}^*_{\rm A} , {\cal L}^*_{\rm B} \right] \right] + \left[{\cal L}^*_{\rm A}, \left[ {\cal L}^*_{\rm O} , {\cal L}^*_{\rm B} \right] \right]  - \tfrac{1}{2} \left[{\cal L}^*_{\rm A}, \left[ {\cal L}^*_{\rm A}, {\cal L}^*_{\rm O} \right] \right]  \right. \\
&& \left. - \tfrac{1}{2} \left[{\cal L}^*_{\rm B}, \left[ {\cal L}^*_{\rm B}, {\cal L}^*_{\rm O} \right] \right] - \tfrac{1}{2} \left[{\cal L}^*_{\rm B}, \left[ {\cal L}^*_{\rm B}, {\cal L}^*_{\rm A} \right] \right]  \right) + { O}(\delta t^4),
\end{eqnarray*}
where recall that
\[
{\cal L}^*_{\rm Det} = {\cal L}^*_{\rm A} + {\cal L}^*_{\rm B},
\]
and
\[
{\cal L}^*_{\rm LD} = {\cal L}^*_{\rm Det} + {\cal L}^*_{\rm O}.
\]
The calculation of the inhomogeneity in (8) then amounts to a straightforward computation of the commutator series applied to $\rho_\beta$. The commutators needed are:
\begin{align*}
\left[{\cal L}^*_{\rm A}, \left[ {\cal L}^*_{\rm A}, {\cal L}^*_{\rm B} \right] \right] \rho_\beta =& 2 \beta p^T U''(x) \nabla U(x) \rho_\beta - \rho_\beta \beta \, p \cdot \nabla_x p^T U''(x) p, \\
\left[{\cal L}^*_{\rm B}, \left[ {\cal L}^*_{\rm B}, {\cal L}^*_{\rm A} \right] \right] \rho_\beta =& -2 \beta p^T U''(x) \nabla U(x) \rho_\beta , \\
\left[{\cal L}^*_{\rm O}, \left[ {\cal L}^*_{\rm A}, {\cal L}^*_{\rm B} \right] \right] \rho_\beta =& 2 \gamma \left(\Delta_x U(x) - \beta p^T U''(x) p \right) \rho_\beta, \\
\left[{\cal L}^*_{\rm A}, \left[ {\cal L}^*_{\rm O}, {\cal L}^*_{\rm B} \right] \right] \rho_\beta =& \gamma \beta \left(  |\nabla U(x) |^2 -  p^T U''(x) p \right) \rho_\beta, \\
\left[{\cal L}^*_{\rm A}, \left[ {\cal L}^*_{\rm A}, {\cal L}^*_{\rm O} \right] \right] \rho_\beta =& 2 \gamma \left( \beta |\nabla U(x) |^2 - \Delta_x U(x) \right) \rho_\beta, \\
\left[{\cal L}^*_{\rm B}, \left[ {\cal L}^*_{\rm B}, {\cal L}^*_{\rm O} \right] \right] \rho_\beta =& 0 \\
\left[{\cal L}^*_{\rm O}, \left[ {\cal L}^*_{\rm O}, {\cal L}^*_{\rm Det} \right] \right] \rho_\beta =& 0,
\end{align*}
where we have abbreviated the Hessian $\nabla_x \nabla_x^T U(x) =: U''(x).$   

Hence we  see directly that
\[
12 {\cal L}^*_2 \rho_\beta = 3 \gamma \left( \Delta_x U(x) - \beta p^T U''(x) p \right) \rho_\beta + 3 \beta p^T U''(x) \nabla U(x) \rho_\beta - \rho_\beta \, \beta p \cdot \nabla_x p^T U''(x) p,
\]
giving
\begin{eqnarray}
{\cal L}^*_2 \rho_\beta &=& \rho_\beta \left[ \frac{\gamma}{4}  \left( \Delta_x U(x) - \beta p^T U''(x) p \right) + \frac{\beta}{4}  p^T U''(x) \nabla U(x)  - \frac{\beta}{12}  \, p \cdot \nabla_x p^T U''(x) p \right]. \label{fredalt}
\end{eqnarray}
Observe that Eq. \ref{FA} is satisfied since the average of the first term is equivalent to a canonical average which vanishes
\[
{ \mathbb{E} \left( \Delta_x U(x) - \beta p^T U''(x) p  \right) = 0},
\]
whereas the other terms in Eq. \ref{fredalt}, being canonical averages of terms which are odd-order in $p$, necessarily also average to zero.

An analogous computation can be performed for the ABOBA method, giving a slightly (but crucially) different perturbation operator ${\cal L}^{* \, ({\rm ABOBA})}$, where
\[
{\cal L}^{* \, ({\rm ABOBA})}_2 \rho_\beta = - \rho_\beta \left[ \frac{1}{4} \gamma \left( \Delta_x U(x) - \beta p^T U''(x) p \right) + \frac{1}{4} \beta p^T U''(x) \nabla U(x)-  \frac{1}{6} \beta\, p \cdot \nabla_x p^T U''(x) p \right]
\]
Here $U''$ is the Hessian matrix of $U$ and $\Delta_x = \sum_{i=1}^N \frac{\partial^2}{\partial x_i^2}$ is the partial Laplacian in $x$.   Equation \ref{FA} is again seen to be satisfied.

\subsection{Expansion in powers of $\gamma^{-1}$}
 Although we are not able to give the general analytical solution to the partial differential equation of Eq.\ \ref{eqn::pde}, we can find a solution in an important limiting case: the high friction regime.   To do this, we expand the invariant density of the numerical method further, viewing both $\delta t$ and $\varepsilon=\gamma^{-1}$ as small parameters:
\begin{equation}
\hat{\rho} = \exp(-\beta [ H + \delta t^2 ( f_{2,0} +  f_{2,1}\varepsilon +O(\varepsilon^2) ) +O(\delta t^4) ]). \label{densityexpansion}
\end{equation}
Dividing by $\gamma$, we may reduce Eq. \ref{eqn::pde} to
\[
[{\cal L}^*_0 + \varepsilon {\cal L}^*_1](f_{2,0} + \varepsilon f_{2,1} + O(\varepsilon^2 )) = g_{0} + \varepsilon g_{1},
\]
where
\[
{\cal L}^*_{0} = \frac{1}{\beta} \Delta_p  - p \cdot \nabla_p, \quad {\cal L}^*_{1} = \nabla U(q) \cdot \nabla_p - p \cdot \nabla_x,
\]
and, for BAOAB,
\begin{align*}
g_0 & = \tfrac{1}{4} \left( \beta^{-1} \Delta_x U(x) - p^T U''(x) p \right),\\
g_1 & = \tfrac{1}{4}  p^T U''(x) \nabla U(x) - \tfrac{1}{12} p \cdot \nabla_x p^T U''(x) p.
\end{align*}
Note that this is a singularly perturbed system as ${\cal L}_0$ is degenerate and it is only the combined operator that has the necessary regularity to define a unique solution.  Nonetheless, as explained below it is possible to find the leading term $f_{2,0}$ by substituting a truncated expansion of fixed degree and solving the resulting equations.     

Equating powers of $\varepsilon$, we find
\[
{\cal L}^*_0 f_{2,0} = g_0, \quad {\cal L}^*_1 f_{2,0} + {\cal L}^*_0 f_{2,1} = g_1,
\]
and
\[
\qquad {\cal L}^*_{1} f_{2,n-1} = - {\cal L}^*_{0} f_{2,n} \, \, \, \textrm{(for $n>1$)}.
\]

Truncating at $n=2$, for example, we find the following solution of these equations:
\begin{align*}
f_{2,0} \equiv f_{2,0}^{\rm BAOAB} & = \tfrac{1}{8}  \left( p^T U''(x) p -  \beta^{-1}\Delta U(x) \right), \\
f_{2,1} \equiv f_{2,1}^{\rm BAOAB} & = \tfrac{1}{24} \beta^{-1} p^T \nabla_x \Delta_x U(x) - \tfrac{1}{72}  p^T \nabla_x p^T U''(x) p,\\
f_{2,2} \equiv f_{2,2}^{\rm BAOAB} &= \tfrac{1}{296} p^T \nabla_x  p^T \nabla_x p^T U''(x) p -\tfrac{1}{48}  \nabla U(x) \cdot  \nabla_x p^T U''(x) p.
\end{align*}

For ABOBA,  $f_{2,0}$ the solution would change to
\[
f_{2,0}^{\rm ABOBA}  = - \tfrac{1}{8}  \left( p^T U''(x) p - 2 \beta^{-1}\Delta U(x) \right).
\]

\subsection{Marginal distribution} 
We now turn out attention to the configurational marginal distribution obtained by integrating the density expansion Eq. \ref{densityexpansion} with respect to the momenta.
Our interest is only in the leading term of this expansion, which defines the sampling behavior for large $\gamma$. 
Ignoring higher order terms in $\delta t$ and $\varepsilon$, we would have the distribution
\[
\hat{\rho} = \tilde{\rho} \times  (1 +{\cal R}),
\]
where, for BAOAB,  $\tilde{\rho} =\rho_{\beta} \times \exp \left ( -\beta \left [\tfrac{\delta t^2}{8}  \left( p^T U''(x) p -  \beta^{-1}\Delta U(x) \right)  \right ]  \right )$,
and ${\cal R} = O(\varepsilon \delta t^2) + O(\delta t^4)$.  
 
Noting that the leading term in the exponent of $\tilde{\rho}$ is quadratic in momenta, we integrate out with respect to $p$ to obtain
\begin{align*}
\int \tilde{\rho} \, {\rm d}^N p & = \int  \exp \left ( -\beta \left [ \frac{1}{2} p^T(I + \frac{\delta t^2}{4}U'') p + U - \beta^{-1} \frac{\delta t^2}{8} \Delta U \right ] \right )  {\rm d}^N p\\
& = \sqrt{2 \pi k_B T /\det\left (I + \frac{\delta t^2}{4}U'' \right ) } \exp \left ( -\beta \left [  U - \beta^{-1} \frac{\delta t^2}{8} \Delta U \right ] \right )
\end{align*}
Using the identity $\det(M) = \exp({\rm trace}(\log(M)))$, we find
\[
\int \tilde{\rho} \, {\rm d}^N p \propto \exp \left ( -\frac{1}{2} {\rm trace} \left ( \log \left ( I + \frac{\delta t^2}{4}U''  \right ) \right ) \right )  \times  \exp \left ( -\beta \left [  U - \beta^{-1} \frac{\delta t^2}{8} \Delta U \right ] \right ) .
\]
We then Taylor expand the logarithm of $I+\frac{\delta t^2}{4}U''$ and take the trace to obtain a cancellation of the $\delta t^2$ terms, giving
\[
\int \tilde{\rho} \, {\rm d}^N p \propto \exp(-\beta U + O(\delta t^4)).
\]
The contribution to the configurational distribution error  due to $\tilde{\rho}$ is $O(\delta t^4)$.  This means that the overall error in the marginal distribution of $\hat{\rho}$ (which includes the neglected factor $1+{\cal R}$) will be $O(\varepsilon \delta t^2) + O(\delta t^4)$.  

If $\varepsilon$ is small (or $\delta t$ is relatively large), the error will be dominated by the quartic term in $\delta t$ and we will observe 4th order accuracy in configurational averages.  For small $\delta t$ the method is always eventually second order.

In the case of the ABOBA method, the remarkable cancellation of the second order errors does not occur and the method always exhibits 2nd order configuration distribution error.

\section{The Limit Method}
We now consider the limit $\gamma\rightarrow \infty$, where the exact solution of  the vector OU process reduces to $p= \sqrt{k_BT} M^{1/2} R$, where $R$ is a vector whose components have a standard normal distribution (Gaussian white noise). {Alternatively, we could consider the limit of the particle mass going to $0$, although this requires a reformulation of Eq. \ref{lang} so that the friction is proportional to the velocity instead of the momentum \cite{FreeEnergy}. Whichever limit is taken, we would expect the ultimate result to be the same.} (Here we have reintroduced the masses in order to present the method, since they may be useful scaling parameters in simulation.)   In the configurations it is straightforward to show that the BAOAB method therefore becomes
\begin{align*}
x_{n+1} & = x_n +\frac{ \delta t^2}{2} M^{-1} F(x_n) + \frac{\delta t}{2} M^{-1} (p_n + p_{n+1})\\
&  = x_n + \frac{\delta t^2}{2} M^{-1} F(x_n) + \frac{\delta t}{2} \sqrt{k_BT} M^{-1/2} (R_n + R_{n+1}),
\end{align*}
where the $R_n$ are vectors of i.i.d.\ standard normal random variables.
Replacing $\delta t^2/2$ by $h$ we arrive at Eq.\ \ref{BD-mod}.  Since the Langevin scheme gives 4th order accurate configurational averages in this limit, we expect the method of Eq.\ \ref{BD-mod} to be second order accurate in modified timestep $h$.  Moreover,  since we completely remove the second order term in the Langevin dynamics configurational density expansion,  we expect {to observe this behaviour across all values of $h$.} 

By contrast the ABOBA scheme gives a much more complicated limit method as $\gamma\rightarrow \infty$ which is not in one-step form.

In the Euler-Maruyama method, the  random perturbations introduced at each step are independent.  In the method of Eq. 2, the random perturbation is a scaling of $Z_n =(R_n+R_{n+1})/\sqrt{2}$; the components $Z_n^{(i)}$ of these are independent of each other and decay linearly with timestep:
\begin{align*}
\langle Z_n^{(i)}, Z_n^{(i)} \rangle & = \frac{1}{2} (\langle R_{n+1}^{(i)},R_{n+1}^{(i)} \rangle + \langle R_{n}^{(i)},R_{n}^{(i)} \rangle ) = 1\\
\langle Z_n^{(i)}, Z_{n-1}^{(i)} \rangle & = \langle R_{n}^{(i)},R_{n}^{(i)} \rangle = 1/2\\
\langle Z_n^{(i)}, Z_{n-k}^{(i)} \rangle & =0, \hspace{0.2in} k=2, 3,\ldots
\end{align*}
Thus, in this new method, we use a colored noise which has characteristics that directly depend on the stepsize, although the noise decorrelates in just a couple of timesteps. {This is therefore no longer a Markov process, however it can be reformulated as such if one considers the appropriate extended space (eg. $y_n=\left[x_n, R_n,R_{n+1} \right]$).} 

\section{Numerical Experiments} 
\begin{figure}[bt]
\begin{center}
\includegraphics[width=3in]{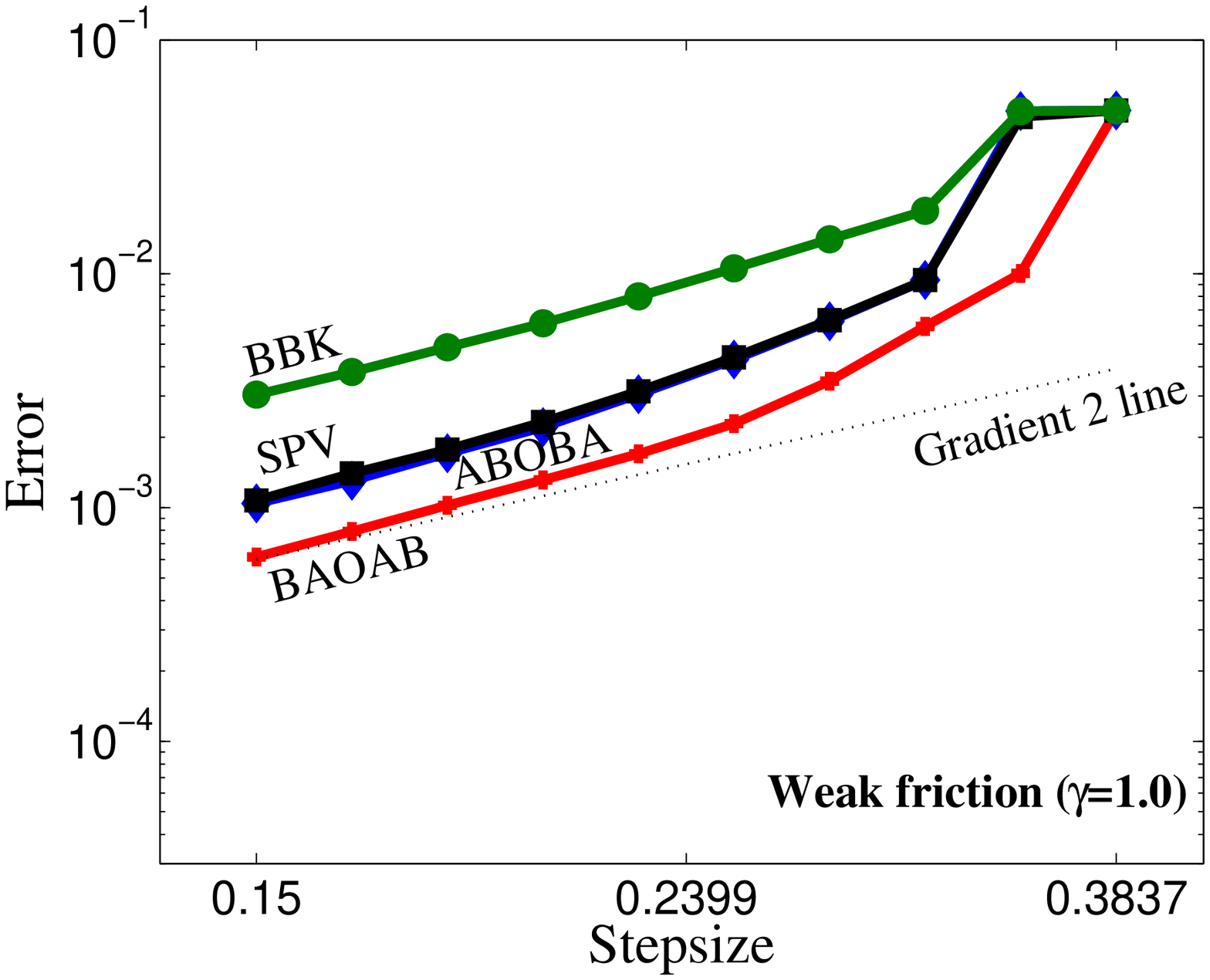}\hspace{0.1in} \includegraphics[width=3in]{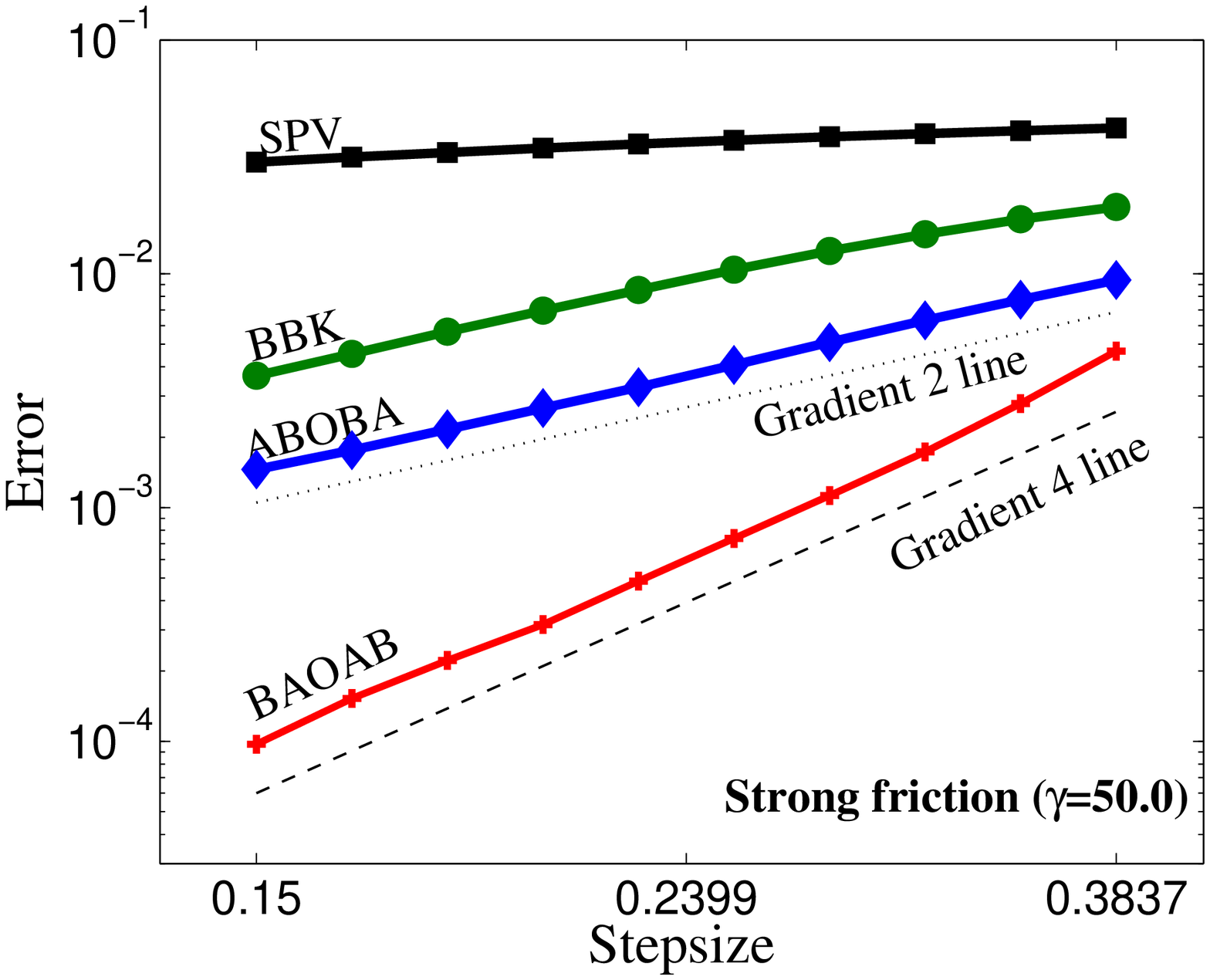}
\caption{The graphs show the comparison of four different Langevin dynamics methods when applied using different stepsizes.  The configurational distribution errors are plotted against the stepsize in a log-log scale. Here $k_BT=1$.  The simulation time was fixed for all runs at $t= 5 \times10^7$, and five runs were averaged to further reduce sampling errors.   At left, $\gamma=1$, at right $\gamma=50$. The graphs are entirely in keeping with the theory presented in the article. \label{quartic}}
\end{center}
\end{figure}

We implemented the methods ABOBA, BAOAB, SPV and BBK and compared the accuracy of configurational sampling for different values of $\gamma$ and a range of timesteps.  A brief analysis shows that the use of the harmonic oscillator leads to special cancellations in the BCH series of the splitting schemes, making it a poor test subject.  Hence, in order to compare the order of accuracy of the different schemes, we first considered an oscillator model in 1D with potential $U(x) = x^4/4 + \sin(1+5x)$. {This was accomplished by introducing $M$ intervals (`bins') of equal length, and computing the mean error in the observed probability frequency compared to the exact expected frequency (obtained by integration of the probability density). If the observed frequency in bin $i$ is ${\omega}_i$, and the exact expected frequency is $\hat{\omega}_i$, then the error calculated is}
\begin{equation}
\rm{Error } = \frac{1}{M} \sum_{i=1}^{M} \left| \omega_i - \hat{\omega}_i \right|. \label{errorformula}
\end{equation}
In this one-dimensional example, we used $20$ bins to cover the interval from $-3.5$ to $3.5$. The configurational density error is plotted against stepsize in log-log scale.  If ${\rm Error} \propto \delta t^r$ then we expect this graph to be a line of slope $r$. {Due to the relative simplicity of this model, we were able to perform highly resolved simulations to calculate accurate error estimates for the configurational distribution. The exact expected value that we compare experimental results against can be computed to arbitrary prescision, and we are able to run as many simulations as needed in order to drive the variance of results to a minimum. The variance in our results, in cases where the stepsize was less than $0.3$, was consistently below $10^{-10}$. Above a stepsize of $0.3$ some of the methods were found to be unstable.}   The results of our simulations are summarized in Figure \ref{quartic}.

\begin{figure}[bt]
\begin{center}
\includegraphics[width=6in]{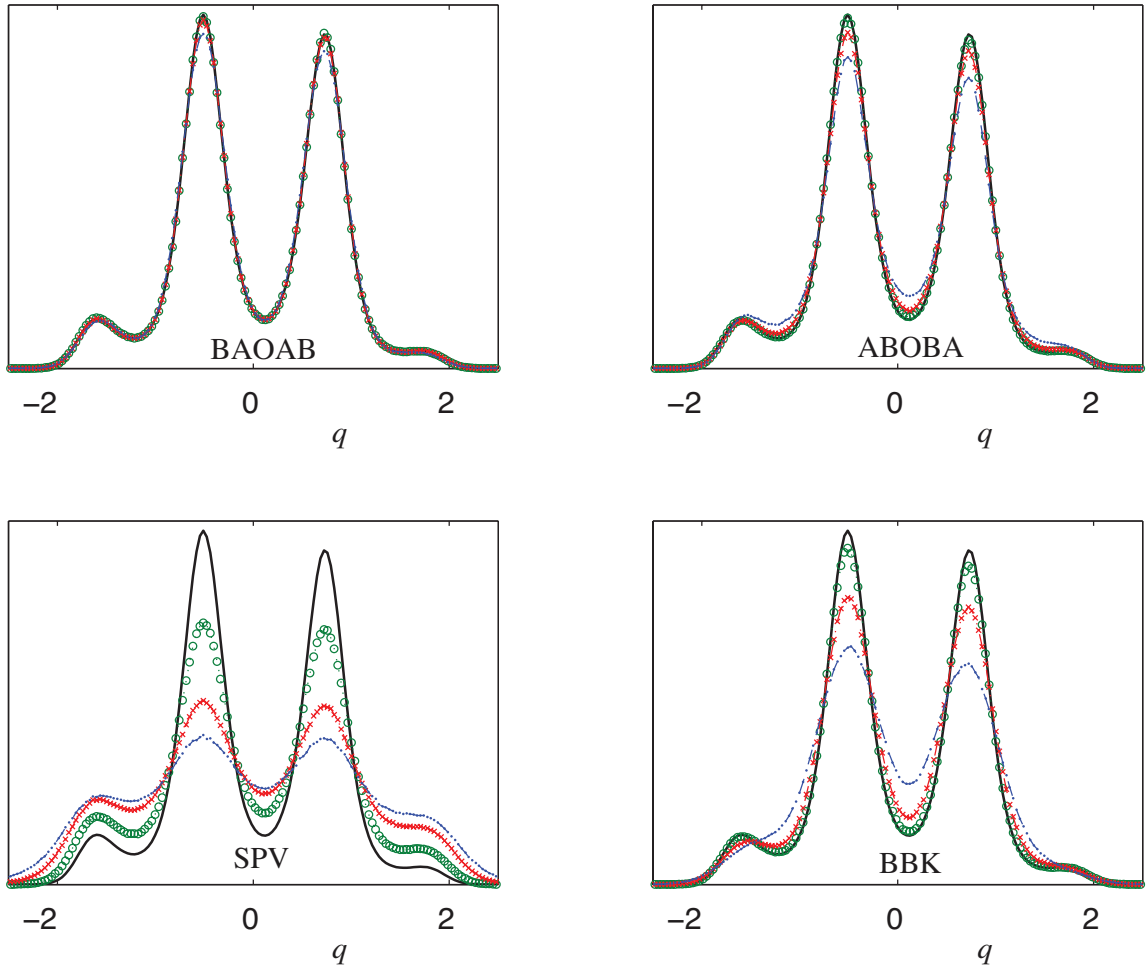}
\caption{Computed distributions of the 1D model problem are compared for $\gamma=20$. The three curves for each of the four methods show the results for three different stepsizes:  $\delta t=0.1$ (circles), $\delta t=0.2$ (crosses), $\delta t=0.3$ (dashed), compared to the dark, solid curve representing the exact distribution.  The graph shows that the generalized GLA methods are superior to SPV and BBK in the moderate $\gamma$ regime.  \label{1ddist}}
\end{center}
\end{figure}

As we can see, when $\gamma$ is small, the methods perform somewhat similarly, at least in the qualitative sense, with all showing a 2nd order error in configurational sampling, and the ABOBA and SPV methods essentially identical.  As $\gamma$ is increased, the {substantial} difference between BAOAB and the other methods becomes apparent.   In the limit of large $\gamma$ the SPV method effectively annihilates the force which results in poor sampling.  In the graph for $\gamma=1$, we can see that for larger values of the stepsize, the graph steepens (indicating that the fourth order term is dominant); as the stepsize is decreased the method exhibits a 2nd order asymptotic decay.   With $\gamma=50$ the fourth order behavior is seen for for all indicated data points, although, again, this becomes second order for smaller values of $\delta t$.    Note that the limit method Eq.\ 2 (with the substitution $h=\delta t^2/2$) gives an essentially identical behavior to the $\gamma=50$ case. We also give a comparison of the actual computed configurational distributions, at different stepsizes using each method, for $\gamma=20$ in Figure \ref{1ddist}.  

\begin{figure}[bt]
\begin{center}
\includegraphics[width=3in]{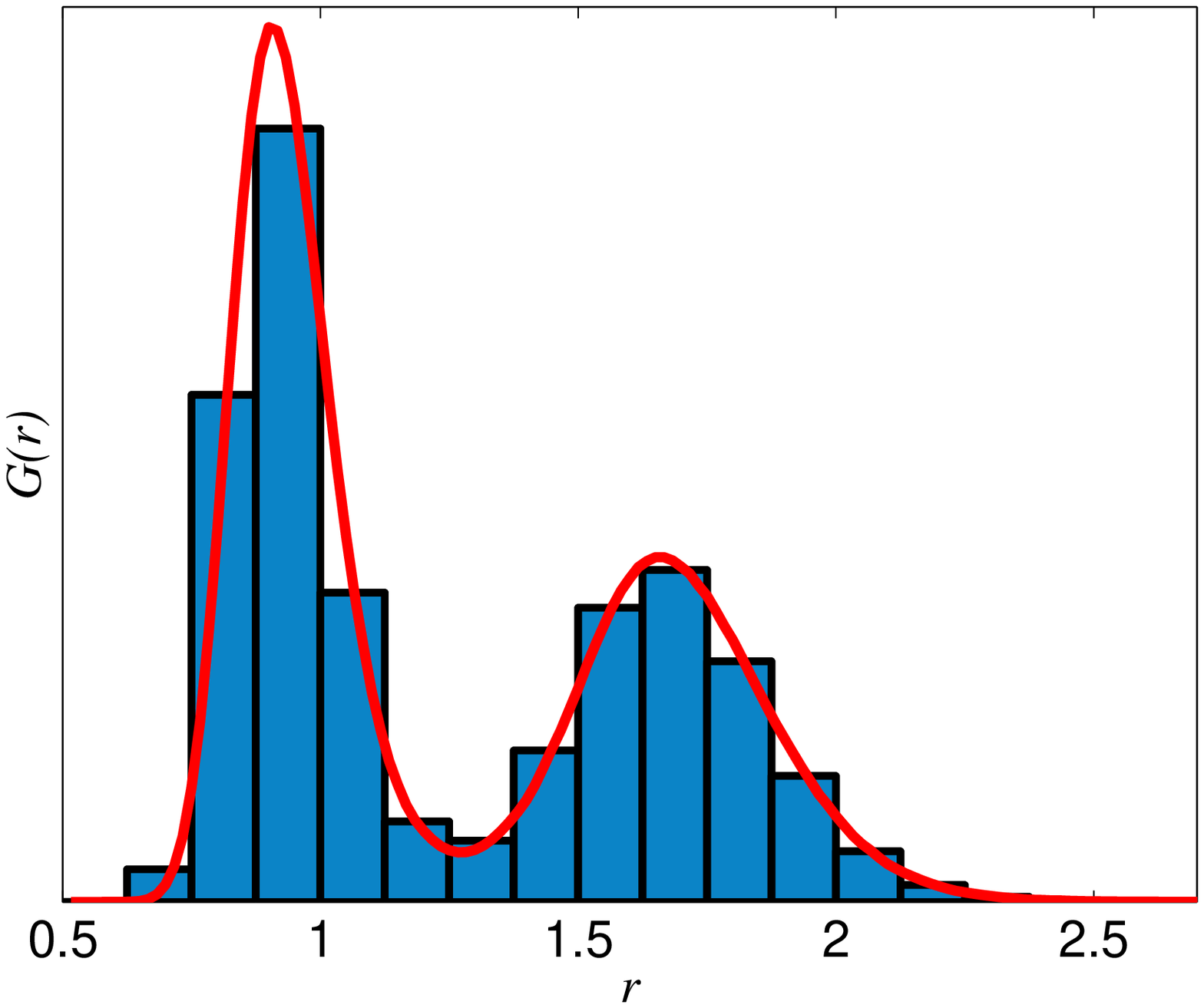}\hspace{0.2in}
\includegraphics[width=3in]{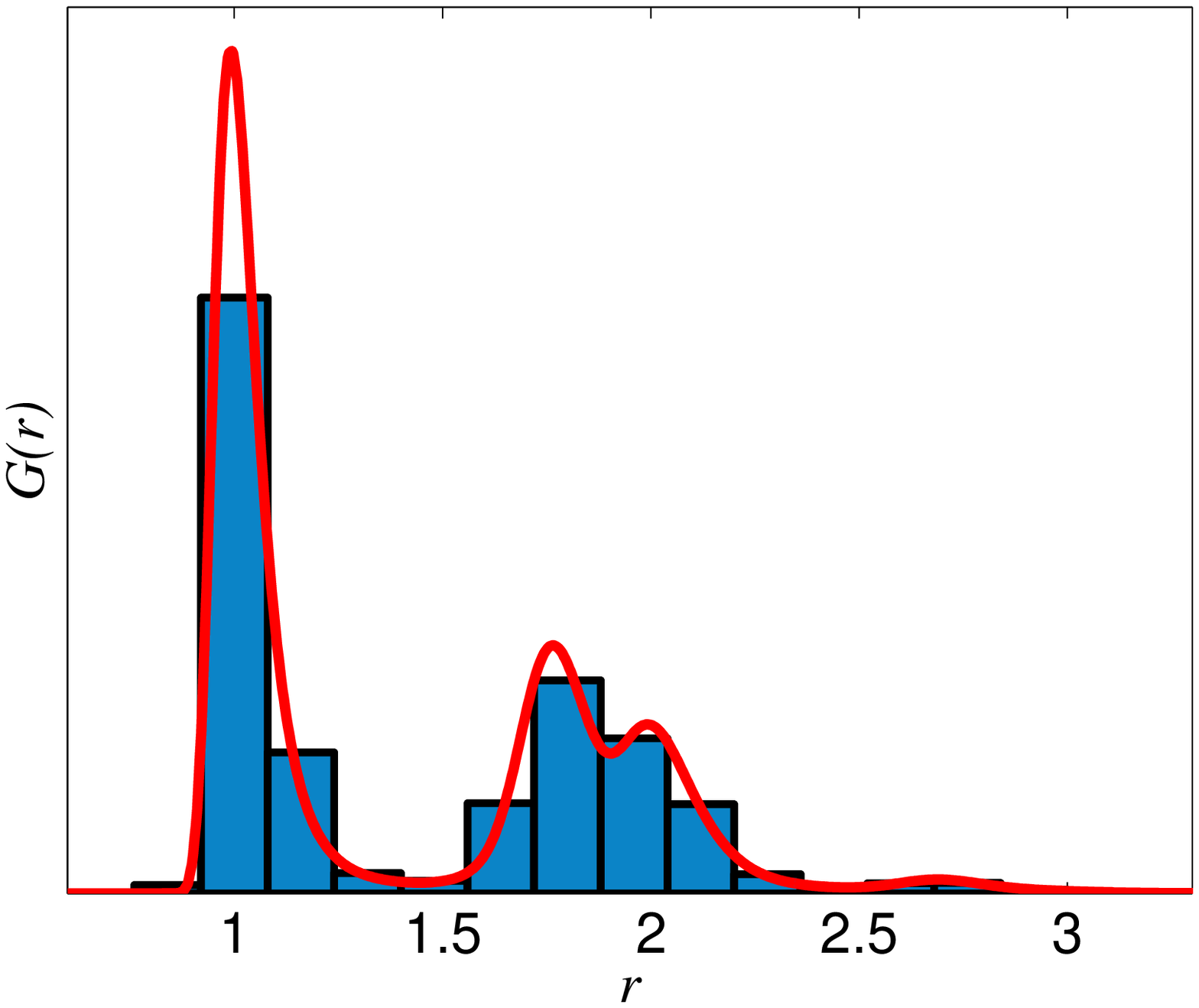}
\caption{The diagrams illustrate the distributions of interatomic distances, $G(r)$, for Morse (left) and Lennard-Jones (right) clusters.     They also show the choice of bins used in calculating the numerical distributions. \label{pic::GofR}}
\end{center}
\end{figure}

To examine the performance of the limit method in more detail, we next considered small molecular clusters consisting of seven atoms (motion restricted to the plane),  with both Morse ($\varphi_M$) and Lennard-Jones ($\varphi_{LJ}$) potentials, given by 
\begin{eqnarray*}
\varphi_{\rm M}(r) &=& (1- \exp(-a(r-r_m)) )^2,\\
\varphi_{\rm LJ}(r) &=& \varepsilon \left( \left( \frac{r_m}{r} \right)^{12} - 2 \left( \frac{r_m}{r} \right)^6 \right),
\end{eqnarray*} 
where we use $a=2, \varepsilon=1$ and $r_m = 1$. The overall potentials are hence
\begin{eqnarray*}
U_{\rm M}(q) &=& \sum_{i=1}^7 \sum_{j>i}^7 \varphi_{\rm M}(r_{ij}),\\
U_{\rm LJ}(q) &=& \sum_{i=1}^7 \sum_{j>i}^7 \varphi_{\rm LJ}(r_{ij}) + \sum_{k=1}^7 \frac{r_k^2}{8},
\end{eqnarray*}
where $r_{ij}$ is the distance between particles $i$ and $j$, $r_k$ is the distance between particle $k$ and the origin, and a mild harmonic term is included in the case of the Lennard-Jones system in order to prevent particles being ejected from the cluster.

The Morse potential gives a smoother dynamics compared to Lennard-Jones (the Morse forces were on average three times smaller than those for Lennard-Jones) and allows a more satisfying determination of the error scaling behavior with stepsize.   To quantify the error in configurational sampling,  we calculated the radial density $G(r)$ by binning the instantaneous interatomic distances at each step into 20 compartments and compared to a calculated reference value (Figure \ref{pic::GofR}).  Though not exact, the errors in the reference will be negligible compared with configurational distribution errors at higher stepsizes. For the Morse cluster the reference stepsize we used was $h_{\rm{ref}}=0.001$, whereas for the Lennard-Jones cluster, we used $h_{\rm{ref}}=0.00025,$ both with the same integration time that was used in their respective test runs. In both cases the cluster was initialized with 6 particles placed on a regular hexagon with unit side-length and with the remaining particle in the centre, and initial velocities randomly drawn from the canonical distribution.

Considerable computation is required to achieve the level of accuracy required, due to the dominance of sampling error and the complexity of the system compared to the earlier 1D example. We ran a large number of independent simulations at each stepsize and computed the average radial density plot for both examples, and compared this to the reference result.

\begin{figure}[b]
\begin{center}
\includegraphics[width=3in]{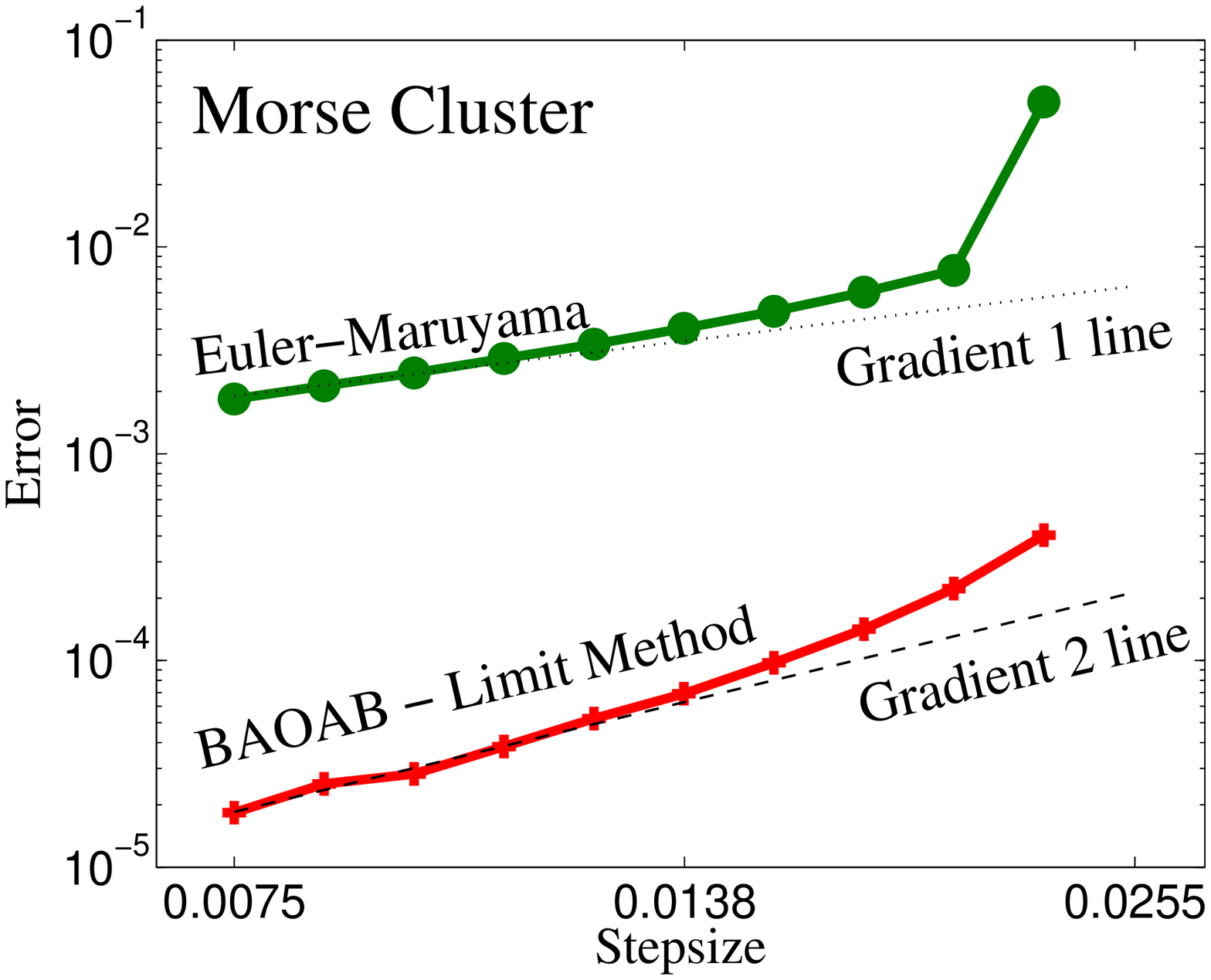} \hspace{0.1in}
\includegraphics[width=3in]{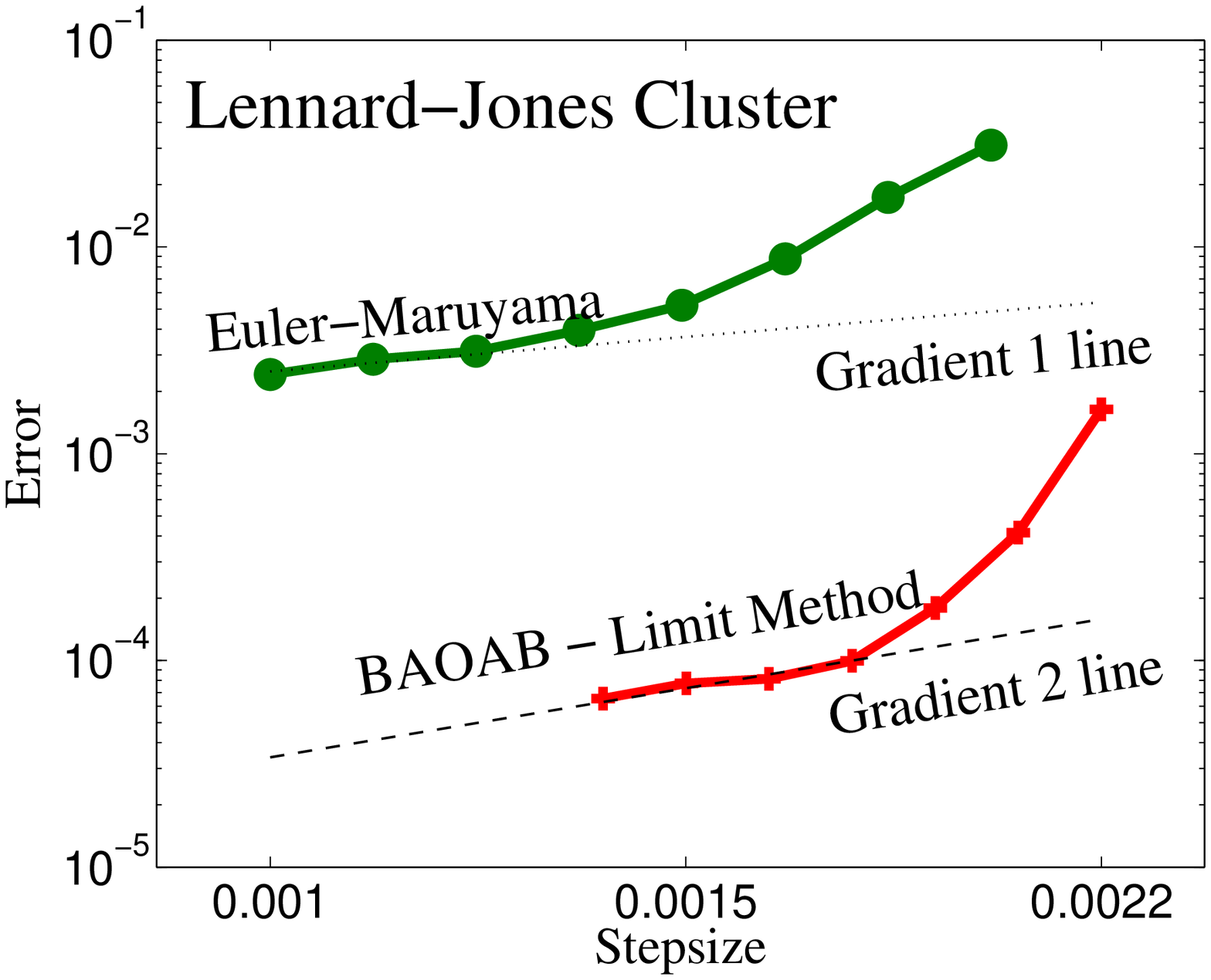}
\caption{The radial distribution errors are plotted in log-log scale against stepsize, demonstrating the first order decay of the error in the case of Euler-Maruyama and the second-order behavior of the BAOAB limit method (in modified timestep $h$).  Left: Morse potential; Right: Lennard-Jones.    For Morse we used a temperature of $k_BT=0.1$, a fixed  time interval of $t=4\times 10^6$, with stepsizes ranging from $0.0075$ to $0.0225$. For Lennard-Jones the temperature was $k_B T=0.2$,  $t=2.5\times 10^5$ and stepsizes ranged from $.001$ to $0.0022$. In order to drive the variance of the results down, a large number of runs were necessary: for the Morse simulation the error is computed using the average histogram computed from 200 independent runs, while the Lennard-Jones simulation used 1000 independent runs.   Around two orders of magnitude of improvement are observed in the accurate regime, but, perhaps even more important, the BAOAB limit method is usable at substantially larger stepsizes than Euler-Maruyama. \label{morse}}
\end{center}
\end{figure}

Our results, presented in Figure \ref{morse}, are entirely consistent with our analysis and show the second order dependence of the configurational sampling error on $h$ (equivalent to fourth order in $\delta t$), as compared to the Euler-Maruyama's first order behavior. These results demonstrate a good agreement with our theoretical results, however even using extensive computation the variances in each experiment were still quite high. If $\omega_{h,n,m}$ is the density of bin $m$ in simulation $n$ at stepsize $h$, we calculate the variance as
\[
 \sigma_h^2 = \frac{1}{N M} \sum_{n=1}^N \sum_{m=1}^M \left( \omega_{h,n,m} - \frac{1}{N} \sum_{k=1}^N \omega_{h,k,m} \right)^2,
\]
where we used $M=20$ bins for each experiment, $N=200$ for the Morse experiment and $N=1000$ for the Lennard-Jones experiment. The variances for the Morse experiment were around $10^{-10}$ while the Lennard Jones experiment had variances around $10^{-8}$. We expect that completing more simulations (increasing $N$) would reduce the variance and give smoother results in Figure \ref{morse}. Nonetheless, in both cases, we observe a significant reduction in the error compared to Euler-Maruyama.


\begin{figure}[ht]
  \centering
  \subfloat[BBK]{\label{fig:rocbbk}\includegraphics[width=3in]{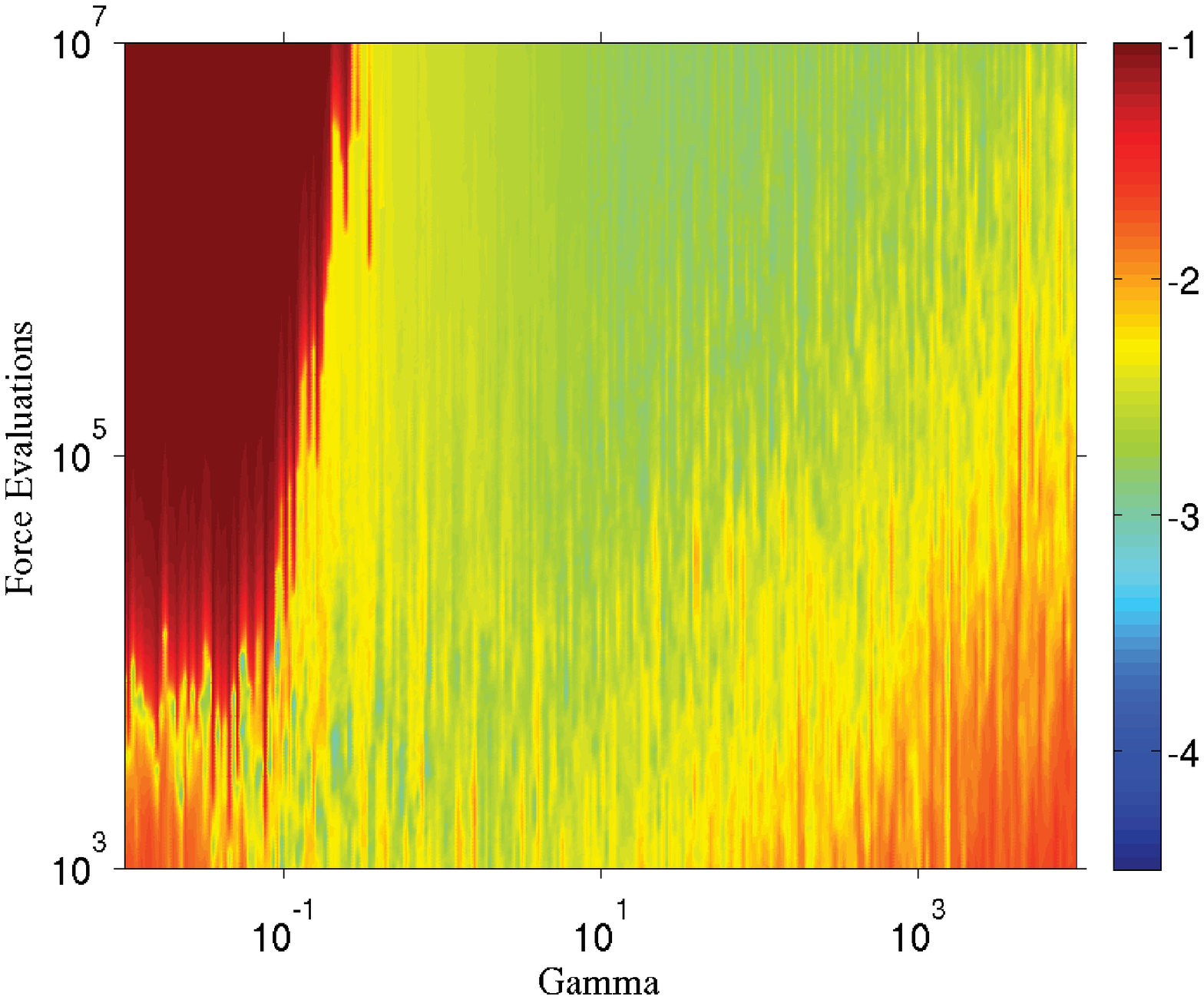}}    \hspace{0.2in}            
  \subfloat[SPV]{\label{fig:rocspv}\includegraphics[width=3in]{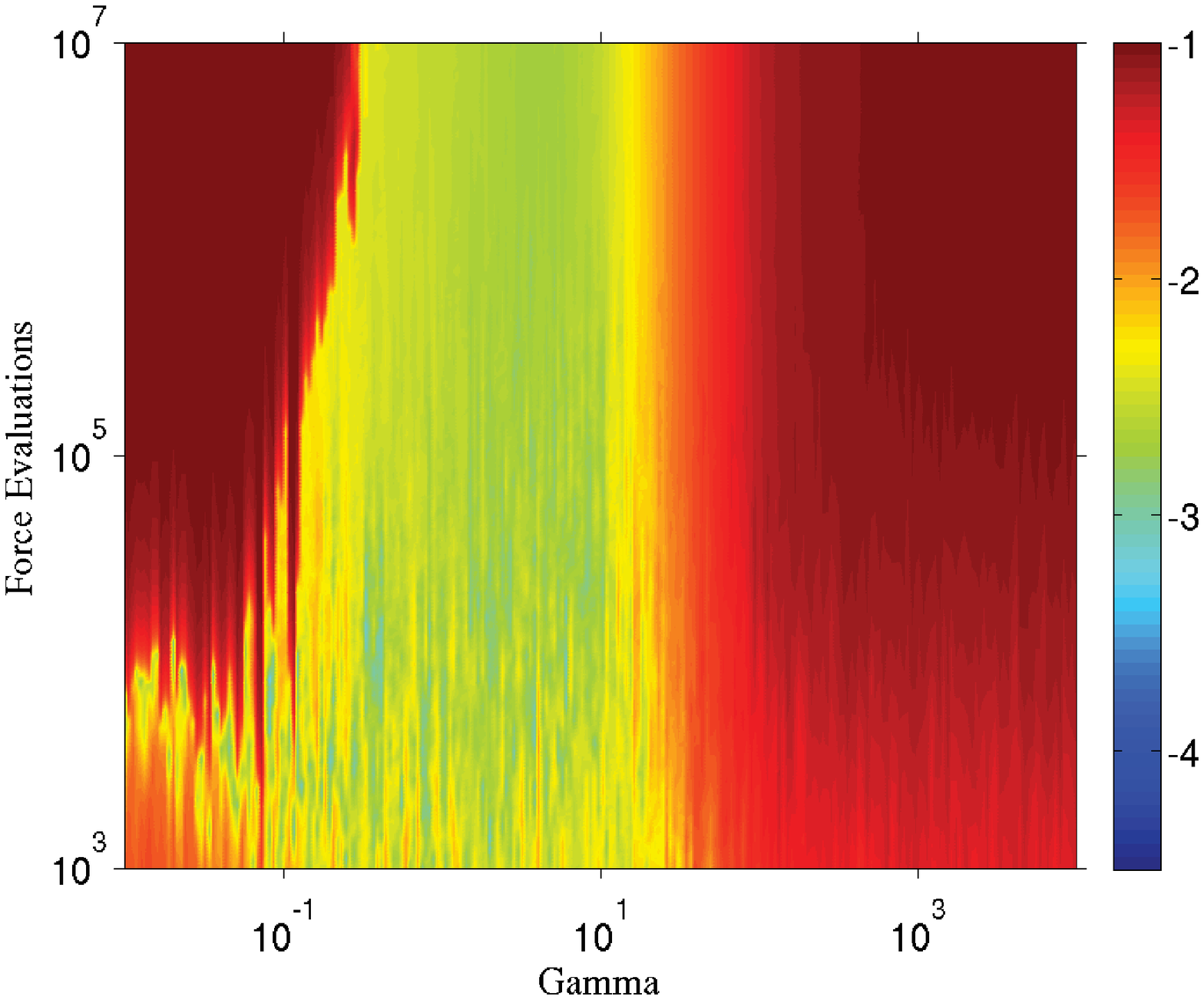}}\\
  \subfloat[ABOBA]{\label{fig:rocaboba}\includegraphics[width=3in]{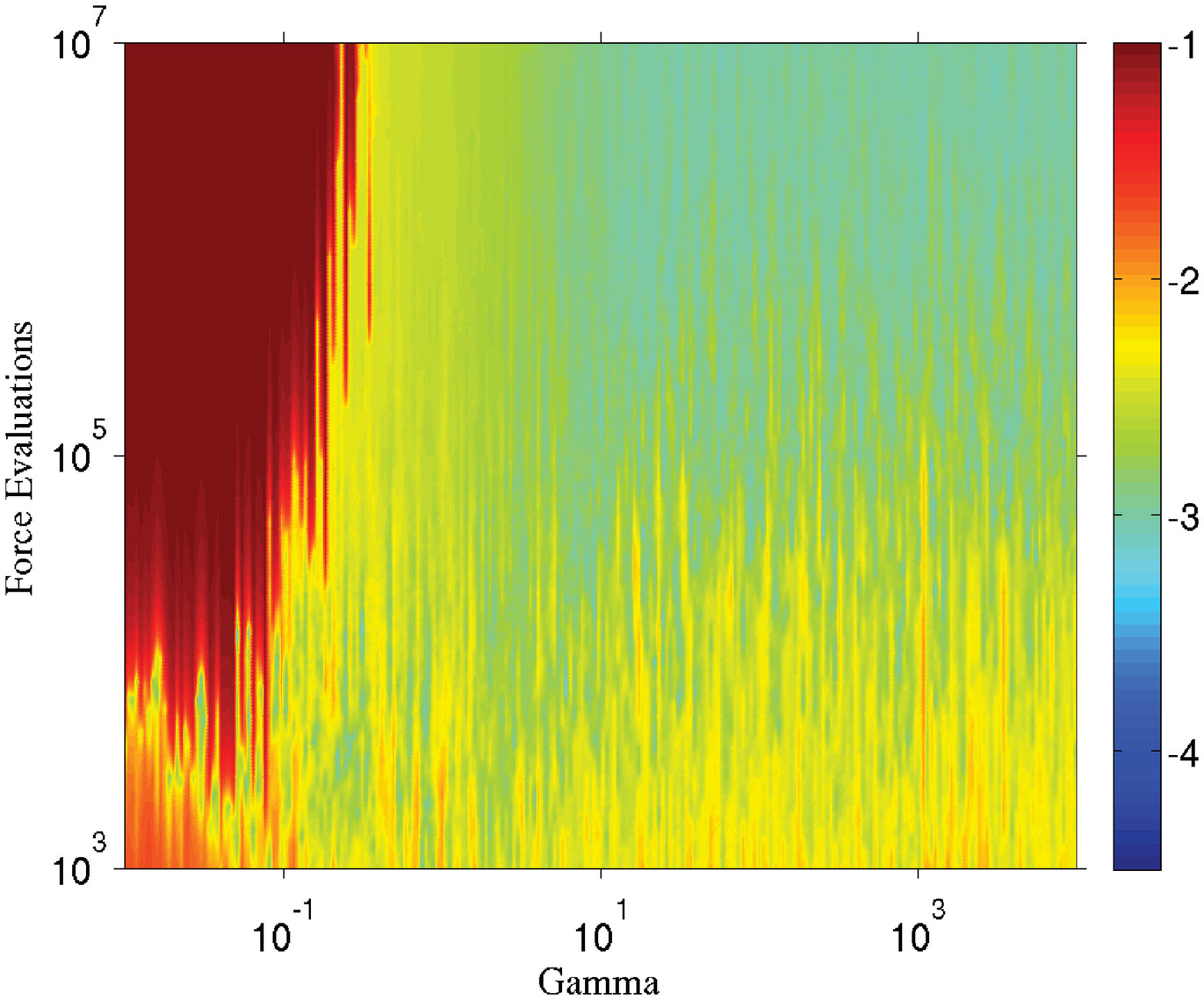}} \hspace{0.2in} 
  \subfloat[BAOAB]{\label{fig:rocbaoab}\includegraphics[width=3in]{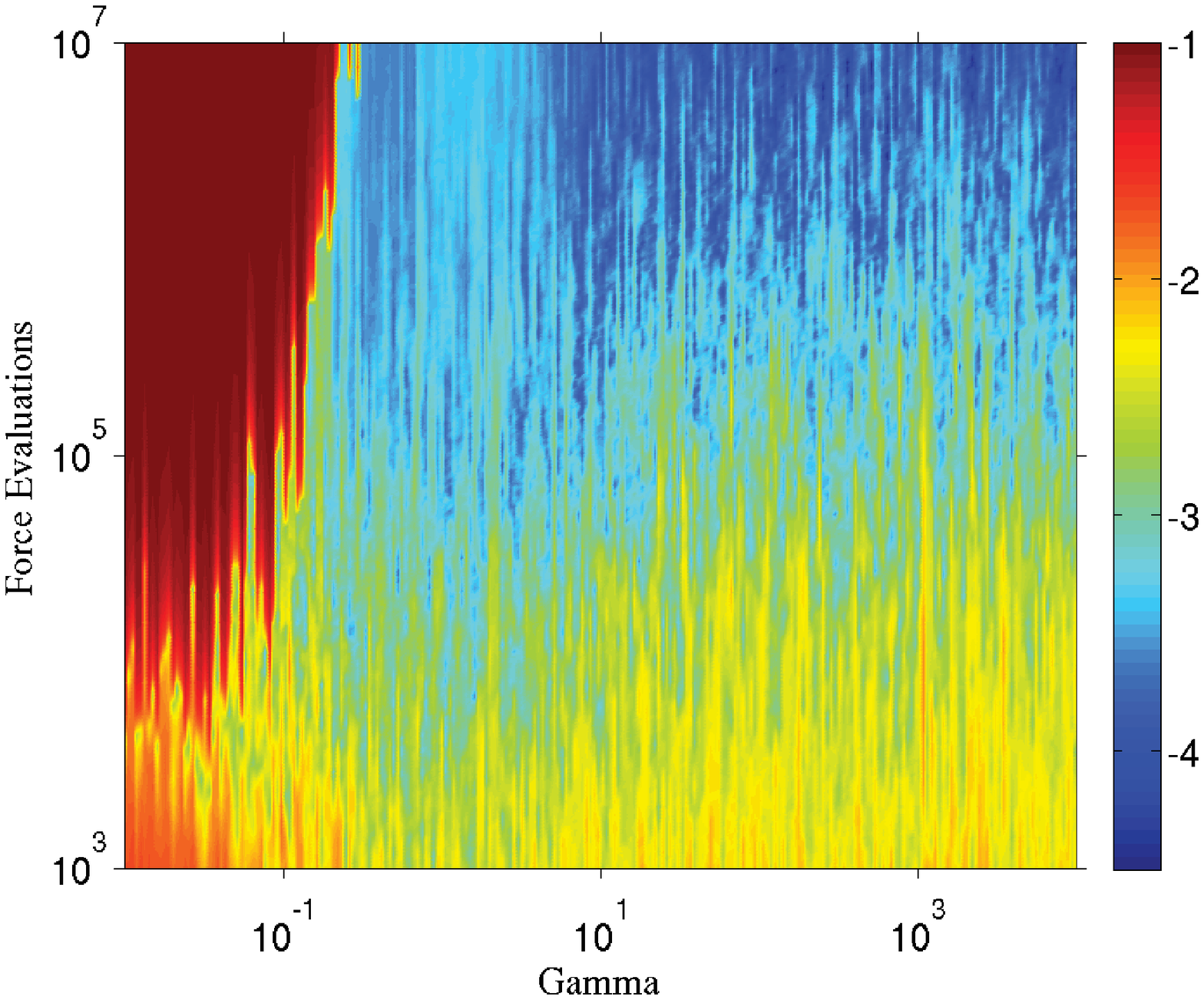}}
\caption{The diagrams illustrate the performance of the different algorithms (labelled) by showing the error in the computed radial distribution functions as a function of both $\gamma$ and the number of timesteps (samples) taken, in the case of the 7 atom Lennard-Jones model.   The graphs address the potential concern that the larger values of $\gamma$ needed to give the superconvergence property may reduce the rate of convergence to equilibrium (it does not, in the case of BAOAB).    The same stepsize of $\delta t=0.044$ was used for all these simulations.} 
  \label{rateofconv}
\end{figure}

It might be a suggested that the improved accuracy seen in the high $\gamma$ regime could potentially come at the price of a slower convergence to equilibrium {due to a reduced rate of transition between metastable states, hence overall sampling of the configurational distribution might be impaired in favour of local sampling.  To address this point, we have plotted in Figure \ref{rateofconv} the error computed in our Lennard-Jones simulation as a function of the number of force evaluations  (vertical) against the friction value  $\gamma$ used for the simulation  (horizontal). Gridpoints in the plot are coloured according to the configurational error, computed using Eq. \ref{errorformula}.   The results indicate that the convergence rate for the BAOAB method is not diminished for large friction coefficient, so it does not appear that we are sacrificing sampling accuracy using this scheme.  Note that the performance of ABOBA is also robust in the limit of large $\gamma$, but the achievable accuracy is reduced as it is only second order, consistent with what we have presented in this article.} 

\section{Conclusions} 
Our results confirm the {theoretical results of Sections 3 and 4, in particular showing the higher order configurational sampling of the BAOAB method}.    The order in its Langevin formulation is effectively four in the large $\gamma$ limit and large values of $\gamma$ do not impair its stability due to the use of exact Ornstein-Uhlenbeck solves.   Not only is the order of the method high, the constant multiplying the leading term must be, in the cases looked at here, of modest magnitude, since the errors are relatively small also at the large stepsize stability threshold.      Since, in the context of molecular dynamics, BAOAB is a `cheap' and easy to implement scheme using only a single force vector per timestep, we stress that there is no price to pay for its improved accuracy.  

In molecular modelling, there are other errors that play important roles, most importantly errors in the force fields (or, more fundamentally, the errors due to not modelling quantum mechanics properly) and sampling errors.  Obviously these errors may dominate the overall method error and limit the relative benefit to be gained by using one integrator as compared to another, but it is also clear that both of the other types of errors are constantly being reduced through the design of better models and the use of more powerful computers.  More important, one can ask the question: how can a practitioner know which part of the error in a given complicated simulation is due to sampling error and which part due to the truncation errors addressed here?    In our experiments with molecular models, even where there was substantial sampling error still present, we nonetheless found the accuracy to be noticeably higher for {the BAOAB method}; it is likely that this improvement in sampling accuracy would be of direct benefit in many real world simulations.   Finally we point out that the BAOAB scheme and its limit method (Eq.\ 2) was, in each case studied, stable at a larger stepsize  than the alternatives, meaning that longer time intervals are made accessible.  This was particularly dramatic in the case of the Lennard-Jones system. 


\vspace{0.2in}
\par\noindent {\bf Acknowledgements.} The first author acknowledges the support of a JTO Fellowship from the Institute for Computational  Engineering and Sciences at the University of Texas. The second author was supported  by the Centre for Numerical Algorithms and Intelligent Software (funded by EPSRC grant EP/G036136/1 and the Scottish Funding Council).   A. Stuart, H. Owhadi, and especially G. Stoltz made helpful suggestions regarding the presentation of the results.  We further acknowledge the comments received from A.  Abdulle and M. Tretyakov which have improved the article and its relation to other works.   Computations were performed on the University of Edinburgh's ECDF compute facility.

\section*{Appendix: Langevin Dynamics Integrators} 
The Langevin dynamics methods used for the numerical experiments in this paper are given here.  {   $R_i$ is an $N-$vector of i.i.d. Normal random numbers,} with mean $0$ and variance $1$. The diagonal mass matrix is denoted $M$, and we assume a timestep $\delta t$ is provided. Given the parameter $\gamma$, we define useful constants
\[
c_1 = e^{-\gamma \delta t}, \quad c_2 = \gamma^{-1} (1-c_1),
\]
and
\[
c_3 = \sqrt{k_BT(1-c_1^2)}.
\]
\par\noindent {\bf `BAOAB' method} 
\begin{flalign*}
p_{n+1/2} &=   p_{n} - \delta t \, \nabla U(x_{n}) /2 ;&\\
x_{n+1/2} &=  x_{n} + \delta t \, M^{-1} p_{n+1/2} /2;&\\
\hat{p}_{n+1/2} &=   c_1  p_{n+1/2} + c_3  M^{1/2} R_{n+1};&\\
x_{n+1} &=   x_{n+1/2} + \delta t \, M^{-1} \hat{p}_{n+1/2} /2;& \\
p_{n+1} &=   \hat{p}_{n+1/2} - \delta t \, \nabla U(x_{n+1}) /2;&
\end{flalign*}

\par\noindent{\bf `ABOBA' method} 
\begin{flalign*}
x_{n+1/2} &=  x_{n} + \delta t \, M^{-1} p_{n} /2;&\\
p_{n+1/2} &=   p_{n} - \delta t \, \nabla U(x_{n+1/2}) /2 ;&\\
\hat{p}_{n+1/2} &=   c_1  p_{n+1/2} + c_3  M^{1/2} R_{n+1};&\\
p_{n+1} &=   \hat{p}_{n+1/2} - \delta t \, \nabla U(x_{n+1/2}) /2;&\\
x_{n+1} &=   x_{n+1/2} + \delta t \, M^{-1} p_{n+1} /2;&
\end{flalign*}

\par\noindent{\bf Stochastic Position Verlet (SPV)} 
\begin{flalign*}
x_{n+1/2} &= x_n + \delta t \, M^{-1} p_n / 2;&\\
p_{n+1} &= c_1 p_n - c_2  \nabla U(x_{n+1/2}) + c_3  M^{1/2} R_{n+1};&\\
x_{n+1} &= x_{n+1/2} + \delta t \, M^{-1} p_{n+1} /2;&
\end{flalign*}

\par\noindent{\bf The Method of Brunger-Brooks-Karplus (1982) (BBK)} 
\begin{flalign*}
p_{n+1/2} &= (1-\delta t \gamma/2)p_n - \delta t \, \nabla U(x_n) /2 + \sqrt{\delta t k_BT  \gamma}M^{1/2} R_n/2;&\\
x_{n+1} &= x_n + \delta t \, M^{-1} p_{n+1/2} ;&\\
p_{n+1} &= [p_{n+1/2} - \delta t \, \nabla U(x_{n+1}) /2 + \sqrt{\delta t k_BT  \gamma} M^{1/2} R_{n+1}/2]/(1+\delta t \gamma/2);&
\end{flalign*}

\par\noindent{ \bf Euler-Maruyama} 
\begin{flalign*}
x_{n+1} &= x_n - \delta t M^{-1}\nabla U(x_n)  + \sqrt{2k_BT \delta t }M^{-1/2}R_n;&
\end{flalign*}

\par\noindent{\bf `BAOAB' - Limit Method} 
\begin{flalign*}
x_{n+1} &= x_n - \delta t M^{-1} \nabla U(x_n) + \sqrt{\frac{k_BT \delta t}{2}} M^{-1/2}(R_n+R_{n+1});&
\end{flalign*}


\begin{thebibliography}{99}
\bibliographystyle{pnas}
\bibitem{BBK}
Br\"unger, A, Brooks~III, C,  \& Karplus, M.
\newblock (1982) Stochastic boundary- conditions for molecular-dynamics simulations of ST2 water. 
{\em Chem. Phys. Letters} {\bf 105}: 495--500.

\bibitem{Ta1995}
Talay, D.\
\newblock (1995)
Simulation and numerical analysis of stochastic differential systems: a review.
in {\em Probabilistic Methods in Applied Physics}, P. Kree and W. Wedig (Eds.), vol. 451 of Lecture Notes in Physics, pp. 63--106. Springer-Verlag.

\bibitem{MiSc96}
Mishra, B. \ and Schlick, S.\
\newblock (1996) The notion of error in Langevin dynamics. I. Linear analysis
{\em J. Chem. Phys.} {\bf 105}:299--318. 

\bibitem{MiTr2004}
Milstein, G.\  and Tretyakov, M.\, 
\newblock (2004) {\em Stochastic numerics for mathematical physics}, {Springer}.

\bibitem{BuLy2009}
Burrage, K.\ and Lythe, G.\
\newblock (2009) Accurate stationary densities with partitioned numerical methods for stochastic differential equations.
{\em SIAM J. Numer. Anal. } {\bf 47}:1601--1618. 

\bibitem{FreeEnergy}
Leli\`{e}vre, T., Rousset, M., \& Stoltz, G.
\newblock (2010) {\em Free energy computations. A mathematical perspective.} Imperial College Press.


\bibitem{ImportanceSampling}
Liu, J.
\newblock (2001) {\em Monte carlo strategies in scientific computing.} Springer Series in Statistics.


\bibitem{LaMaLe2007}
Larini, L. Mannella R., \&  Leporini, D.
\newblock (2007) Langevin stabilization of molecular-dynamics simulations of polymers by means of quasisymplectic algorithms. {\em J. Chem. Phys.} {\bf 126}:104101.

\bibitem{SkIz02}
Skeel, R and Izaguirre, J.
\newblock (2002) An impulse integrator for Langevin dynamics {\em Mol. Phys.} {\bf 100}:3885--3891.

\bibitem{MiTr2003}
Milstein, G.\  and Tretyakov, M.\, 
\newblock (2003) Quasi-symplectic methods for Langevin-type equations {\em IMA J. Num. Anal.} {\bf 23}: 593--626.



\bibitem{Sh2006}
Shardlow, T., 
\newblock (2006) Modified equations for stochastic differential equations. {\em BIT} {\bf 46}, 111--125

\bibitem{Me2007}
Melchionna, S.\
\newblock (2007) Design of quasi-symplectic propagators for Langevin dynamics.  {\em J. Chem. Phys.} {\bf 127}:044108.

\bibitem{BoOw2010}
Bou-Rabee, N. and Owhadi, H.
\newblock (2010) Long-run accuracy of variational integrators in the stochastic context.
{\em SIAM J.  Num. Anal.} {\bf 48}: 278--297.

\bibitem{Da2010}
Davidchack, R.L.\
\newblock (2010) Discretization errors in molecular dynamics simulations with deterministic and stochastic thermostats. 
{\em J. Comput. Phys.} {\bf 229}: 9323 -- 9346.

\bibitem{Zy2011}
Zygalakis, K.
\newblock (2011) On the existence and applications of modified equations for stochastic differential equations. 
{\em SIAM J. Sci. Comput.} {\bf 33}:102--130.


\bibitem{HaLuWa2006}
Hairer, E. Lubich, C.  \& Wanner, G.
\newblock (2006) {\em Geometric numerical integration}. Springer (New York), Second edition.

\bibitem{LeRe2005}
Leimkuhler, B. and Reich, S.
\newblock (2005) {\em Simulating Hamiltonian dynamics},   Cambridge University Press.

\bibitem{FrSm2001}
Frenkel, D. and Smit, B.
\newblock (2001) {\em Understanding Molecular Simulation}, 2nd Edition, Academic Press.

\bibitem{Ho}
H{\"o}rmander, L.
\newblock (1967) Hypoelliptic second order differential equations. {\em Acta Math.} {\bf 119}: 147--171.

\bibitem{MaStHi2002}
Mattingly, J.C.\, Stuart, A.M.\,  \& Higham, D.J.\
\newblock (2002) Ergodicity for {SDE}s and approximations: locally Lipschitz vector fields and degenerate noise. {\em Stoch. Proc. Appls.} {\bf 101}:185--232.

\bibitem{Vil2009}
Villani, C. 
\newblock (2009) Hypocoercivity. {\em Mem. Amer. Math. Soc.} {\bf 202}

\bibitem{PavHai2008}
Pavliotis, G. and Hairer, M. 
\newblock (2008) From ballistic to diffusive behavior in periodic potentials. {\em J. Stat. Phys.} {\bf 131} 175--202


\bibitem{JouSto2012}
Joubaud, R. and Stoltz, G.
\newblock (2012) Nonequilibrium shear viscosity computations with Langevin dynamics. {\em SIAM MMS} {\bf 10}: 191--216


\bibitem{SiChCr}
Sivak, D. Chodera, J. and  Crooks, G.
\newblock { (2011) Driven Langevin dynamics: heat, work and pseudo-work. {\em ArXiV} ?}

\bibitem{Tal2002}
Talay, D.
\newblock (2002) Stochastic Hamiltonian dissipative systems: exponential convergence to the invariant measure, and discretization by the implicit Euler scheme. {\em  Markov Processes and Related Fields} {\bf 8}(2) 163--198.


\bibitem{Mat2010}
Mattingly, J.C.\, Stuart, A.M.\,  \& Tretyakov, M.
\newblock (2010) Convergence of numerical time-averaging and stationary measures via Poisson equations. {\em SIAM J.  Num. Anal.} {\bf 48}: 552--577.

\bibitem{DeFa2011}
Debussche, A. and  Faou, E.
\newblock { (2011) Weak backward error analysis for SDEs. {\em SIAM J.  Num. Anal.} (to appear).}










%
%
%
%
%
%
%









%
%
%
%
%


\end{thebibliography}
\end{document}